\documentclass{amsart}

\usepackage{amsfonts}
\usepackage{amssymb}

\newtheorem{introthm}{Theorem}

\newtheorem{theorem}{Theorem}[section]
\newtheorem{lemma}[theorem]{Lemma}
\newtheorem{proposition}[theorem]{Proposition}
\newtheorem{corollary}[theorem]{Corollary}

\theoremstyle{definition}
\newtheorem{definition}[theorem]{Definition}
\newtheorem{example}[theorem]{Example}
\newtheorem{remark}[theorem]{Remark}

\def\Chi{{\rm Char}}
\def\mal{\! \cdot \!}
\def\rq#1{\widehat{#1}}
\def\b#1{\overline{#1}}
\def\CC{{\mathbb C}}
\def\KK{{\mathbb K}}
\def\TT{{\mathbb T}}
\def\ZZ{{\mathbb Z}}
\def\RR{{\mathbb R}}
\def\SS{{\mathbb S}}
\def\NN{{\mathbb N}}
\def\QQ{{\mathbb Q}}

\def\Pic{{\rm Pic}}
\def\cone{{\rm cone}}
\def\Supp{{\rm Supp}}
\def\Spec{{\rm Spec}}

\title[Equivariant embeddings]%
{Equivariant embeddings \\
into smooth toric varieties}

\author[J.~Hausen]{J\"urgen Hausen} 
\address{Fachbereich Mathematik und Statistik, Universit\"at Konstanz,
  78457 Konstanz, Germany}
\email{Juergen.Hausen@uni-konstanz.de}

\subjclass{14E25, 14C20, 14L30, 14M25}

\begin{document}

\begin{abstract}
We characterize embeddability of algebraic varieties into smooth toric
varieties and prevarieties. Our embedding results hold also in an
equivariant context and thus generalize a well-known embedding
theorem of Sumihiro on quasiprojective $G$-varieties. The main idea is
to reduce the embedding problem to the affine case. This is done by
constructing equivariant affine conoids, a tool which extends the
concept of an equivariant affine cone over a projective $G$-variety to
a more general framework. 
\end{abstract}

\maketitle

\section*{Introduction}

Classical algebraic geometry mainly deals with quasiprojective
varieties. These varieties thus come embedded into an ambient space
whose structure is rather well understood, a fact on which rely 
many basic ideas and explicit working tools of the classical theory.
The modern concept of defining a variety by gluing affine pieces is
much more flexible, but the price for this flexibility is the loss of
structural insight offered by the ambient space. The intention of
embedding theorems is to regain such insight.

Since about 1970, toric varieties have been thoroughly
studied. A most remarkable feature is their explicit description by
combinatorial data. The class of toric varieties contains the
classical ambient spaces, namely the affine and projective spaces, but
it is considerably larger. In fact, by a theorem of W\l odarczyk
\cite{Wl}, toric varieties may serve as ambient spaces for
surprisingly many varieties: A normal variety $X$ can always
be embedded into a toric prevariety, and $X$ admits an embedding into
a separated toric variety if and only if every two points of $X$ have
a common affine neighbourhood.

In the present article, we study some problems arising from  
W\l odarczyk's result. The first one concerns singularities: On the
one hand, one would like to get rid of the assumption of $X$ being
normal, on the other hand it is important to know when one can choose
a smooth ambient space. So it is natural to ask,
compare~\cite[Problems~5.4 and~5.5]{Wl}: Which varieties admit
embeddings into smooth toric varieties? A second point is the problem
of embedding equivariantly with respect to algebraic group
actions. Such embeddings are for example interesting in the context of
quotient constructions, as these are quite well understood in the
toric case.

To address the above problems, we introduce a tool that generalizes
the concept of an affine cone over a projective variety: An {\it affine
conoid} over a not necessarily projective variety $X$ is
an affine variety $\b{X}$ together with an action of an algebraic
torus $H$ and a dense open invariant subset $\rq{X} \subset \b{X}$
where $H$ acts freely with a geometric quotient $q \colon \rq{X} \to X
= \rq{X}/H$. These affine conoids are the key to reduce the embedding
problem to the affine case. However, they might also be of interest
independently from our applications to embeddings, see
e.g. Remark~\ref{GIT}.

Our first main result characterizes existence of affine
conoids for arbitrary varieties and relates it to embeddability. 
Following Borelli~\cite{Bo}, we call an irreducible variety $X$ {\it
divisorial}, if for every $x \in X$ there is an effective Cartier
divisor $D$ on $X$ such that $X \setminus \Supp(D)$ is an affine
neighbourhood of $x$. The class of divisorial varieties considerably
extends the class of quasiprojective varieties; for example it
includes all $\QQ$-factorial varieties and their subvarieties.
We prove (see Theorem~\ref{goal1}): 
 
\begin{introthm}
For an irreducible variety $X$, the following
  statements are equivalent:
\begin{enumerate}
\item $X$ is divisorial.
\item There exists an affine conoid over $X$.
\item $X$ admits a closed embedding into a smooth toric prevariety of
  affine intersection.
\end{enumerate}
\end{introthm}

Here we say that a prevariety $Y$ is of {\it affine
intersection} if for any two open affine subvarieties of $Y$, their
intersection is again affine. This means that the non-separatedness
of $Y$ is of quite mild nature. In fact, using the appropriate
formulation of divisoriality, we obtain the above result even for
reducible varieties $X$. For the crucial part, the implication
``i)$\Rightarrow$ii)'', we extend known constructions by Cox~\cite{Co}
and Kajiwara~\cite{Ka} from the setting of toric varieties to
arbitrary divisorial varieties. As to the question of embeddability
into separated smooth toric varieties, we obtain (see
Corollary~\ref{intosepchar}):

\begin{introthm}
An irreducible variety $X$ admits a closed embedding into a smooth
  toric variety if and only if for any two $x, x' \in X$, there is an
  effective Cartier divisor $D$ on $X$ such that $X \setminus
  \Supp(D)$ is affine and contains $x$ and $x'$.
\end{introthm}

We now turn to the second problem, namely, equivariant embeddings. If a
connected linear algebraic group acts on a normal quasiprojective
variety, then Sumihiro's Equivariant Embedding
Theorem~\cite[Theorem~1]{Su} guarantees existence of a locally closed
equivariant embedding into a projective space. We extend this
result to the divisorial case (see Theorem~\ref{wlodeq} and
Corollary~\ref{V2eq}):

\begin{introthm}
Let $X$ be a normal divisorial variety with a regular action of a
connected linear algebraic group $G$. 
\begin{enumerate}
\item There exist a smooth toric prevariety $Z$ of affine intersection
  with a linear $G$-action and a $G$-equivariant closed embedding $X
  \to Z$.
\item If for any two $x,x' \in X$, there is an effective Cartier
  divisor $D$ on $X$ such that $X \setminus \Supp(D)$ is affine and
  contains $x$ and $x'$, then one can choose $Z$ to be a separated
  smooth toric variety.
\end{enumerate}
\end{introthm}

Similar to the case of $Z$ being a projective space, a
linear action on a toric prevariety $Z$ is an action induced by some
linear representation ``over'' $Z$; for the precise formulation see
Section~\ref{affconoids}. Again, the basic step in the proof is the
construction of affine conoids, but now in an equivariant manner.
A consequence of Theorem~3 is that every $\QQ$-factorial toric variety
can be embedded into a smooth one by means of a toric morphism (see
Corollary~\ref{toricemb}). 

The criterion of Theorem~2 and Theorem~3 ii) on pairs of
points $x,x' \in X$ can also be formulated for $k$-tuples of
points; we call the resulting property {\it $k$-divisoriality}. In
view of the Kleiman-Chevalley-Criterion, increasing $k$ means
``approximating'' quasiprojectivity. We show that Theorem~2 and
Theorem~3 ii) have analogous statements for $k > 2$, that means
$k$-divisorial varieties can be embedded into $k$-divisorial smooth toric
varieties (see Theorems~\ref{kdivemb} and~\ref{kdiveqemb}). In
particular, we prove a conjecture of W\l odarczyk~\cite[5.3]{Wl} in the
$\QQ$-factorial case.

The present paper is organized as follows: In
Section~\ref{affconoids}, we introduce equivariant affine conoids and
show that they give rise to equivariant embeddings into toric
prevarieties. In Section~\ref{amplegroups} we present our construction
of equivariant affine conoids over divisorial $G$-varieties. The first
main results are proved in Section~\ref{charconoids}. Moreover, we
relate embeddings via affine conoids to classical projective
embeddings, and discuss a consequence concerning Geometric Invariant
Theory in this section. Finally, Sections~\ref{akproperty}
and~\ref{intosep} are devoted to the problem of embedding into
separated and, more specially, $k$-divisorial smooth toric varieties.

\section{Affine conoids and embeddability}\label{affconoids}

Affine cones are a useful tool to study projective varieties. The
purpose of this section is to extend that tool to more general
varieties. We introduce the notion of an equivariant affine conoid
over a $G$-variety $X$, and we show in Proposition~\ref{cone2embed}
that such an affine conoid gives rise to a $G$-equivariant embedding
of $X$ into a certain smooth toric prevariety.

Throughout the whole article, we work in the categories of varieties
and prevarieties defined over a fixed algebraically closed field $\KK$. 
For the general background, we refer for example to~\cite[Chapter~I]{Hu}.
We say that a prevariety $X$ is of {\it affine intersection} if the
diagonal morphism $X \to X \times X$ is affine. A variety is a
separated but possibly reducible prevariety.

Let us recall some notation on group actions. 
A $G$-variety is a variety $X$ together with a
regular action $G \times X \to X$ of an algebraic group $G$.
We say that the action of a $G$-variety $X$ is {\it free at $x \in X$}
if the orbit map $g \mapsto g \mal x$ is a locally closed
embedding of $G$ into $X$. Moreover, we call an action {\it  free} if
it is free at every point.

We shall be concerned with the following type of quotients: A {\it
geometric quotient} for a $G$-variety $X$ is an affine regular map
$p \colon X \to Y$ onto a variety $Y=X/G$ such that the $p$-fibres
are precisely the $G$-orbits and $\mathcal{O}_{Y} =
p_{*}(\mathcal{O}_{X})^{G}$ holds. Sometimes we allow in this
setting also non-separated quotient spaces $Y$; then we speak of
geometric {\it prequotients}. 

\begin{definition}\label{conoiddef}
An {\it affine conoid} over a variety $X$ is an affine variety $\b{X}$
together with a regular action of an algebraic torus $H$ and a dense
open $H$-invariant subset $\rq{X} \subset \b{X}$ where $H$ acts freely
with geometric quotient $q \colon \rq{X} \to X = \rq{X}/H$.
\end{definition}

Clearly this concept includes the classical notion of an affine cone
over a projective variety. In order to present non-projective complete
varieties admitting an affine conoid, we consider toric
varieties. Recall at this point that a {\it toric variety} is a
normal, and hence irreducible variety together with a regular action
of an algebraic torus having a dense free orbit. There exist
many non-projective complete smooth toric varieties (see
e.g.~\cite{Fu}, p.~74), and a construction of Cox~\cite{Co} provides
affine conoids in that cases:

\begin{example}
Let $X$ be a complete smooth toric variety arising from a fan
$\Delta$ in a lattice $N$. Denote by $\Delta^{(1)}$ the set of
onedimensional cones of $\Delta$. Consider the lattice homomorphism
$$ Q \colon \ZZ^{\Delta^{(1)}} \to N, \qquad e_{\varrho} \mapsto v_{\varrho},
$$
where $e_{\varrho}$ is the canonical base vector corresponding to
$\varrho \in \Delta^{(1)}$ and $v_{\varrho}$ denotes the primitive
lattice vector of $\varrho \in \Delta^{(1)}$. For every cone $\sigma
\in \Delta$ let $\sigma^{(1)}$ be the set of its extremal rays and
define a cone 
$$ \rq{\sigma} \; := \; \cone(e_{\varrho}; \; \varrho \in \sigma^{(1)})
\; \subset \; \RR^{\Delta^{(1)}}.  $$
These cones form a fan in $\ZZ^{\Delta^{(1)}}$, and the associated
toric variety $\rq{X}$ is an open subvariety of $\b{X} :=
\KK^{\Delta^{(1)}}$. The toric morphism $q \colon \rq{X} \to X$
defined by $Q$ is a geometric quotient for the free action of the
algebraic torus $H := \ker(q)$ on $\rq{X}$.
\end{example}

As we are also interested in the equivariant setting, we have to fix
an appropriate equivariant notion of an affine conoid. Let $G$ be an
algebraic group, and let $X$ be a $G$-variety. Suppose that $\b{X}$ is
an affine conoid over $X$, and let $H$, $\rq{X}$ and $q \colon \rq{X}
\to X$ be the associated data as in~\ref{conoiddef}. Assume moreover
that $G$ acts also regularly on $\b{X}$.

\begin{definition}\label{equivconedef}
We say that $\b{X}$ is a {\it $G$-equivariant} affine conoid over $X$
  if the actions of $G$ and $H$ on $\b{X}$ commute, $G$ leaves $\rq{X}
  \subset \b{X}$ invariant, and the map $q \colon \rq{X} \to X$ is
  $G$-equivariant.
\end{definition}

In the subsequent constructions, we shall use a characterization of
free torus actions in terms of certain regular functions. Assume that
an algebraic torus $H$ acts regularly on a variety $X$. Recall that a
function $f \in \mathcal{O}(X)$ is called {\it homogeneous} with
respect to a character $\chi \in \Chi(H)$ if $f(t \mal x) = \chi(t)
f(x)$ holds for every $t \in H$ and every $x \in X$.

\begin{remark}\label{freetorus}
Let $H$ be an algebraic torus, and let $X$ be an affine
$H$-variety. The action of $H$ is free at $x \in X$ if and only if
$x$ has an $H$-invariant open neighbourhood $U\subset X$ admitting
for every $\chi \in \Chi(H)$ a $\chi$-homogeneous $f \in
\mathcal{O}(U)$ with $f(x) \ne 0$. 
\end{remark}

We begin the construction of equivariant embeddings with two auxiliary
results concerning the following situation: Let $G$ be a linear
algebraic group and let $Y$ denote an affine $G$-variety. Suppose that
$H$ is an algebraic torus contained as a closed subgroup in the center
of $G$, and that $V \subset Y$ is a $G$-invariant open subset such
that $H$ acts freely on $V$. Under these assumptions we have:

\begin{lemma}\label{cutgeomquots}
There exist a linear $G$-action on some $\KK^{n}$, a $G$-equivariant
closed embedding $\Phi \colon Y \to \KK^{n}$ and an open subset $U
\subset \KK^{n}$ with the following properties:
\begin{enumerate}
\item $U$ is invariant under the actions of $G$ and $\TT^{n} :=
  (\KK^*)^{n}$.
\item $H \subset G$ acts diagonally on $\KK^{n}$ and freely on $U$.
\item $V = \Phi^{-1}(U)$ holds.
\end{enumerate}
Moreover, if $G$ is an algebraic torus, then one can achieve that $G$
acts diagonally on $\KK^{n}$.
\end{lemma}

\proof Choose generators $f_{1}, \ldots, f_{r}$ of $\mathcal{O}(Y)$
such that for some $s<r$, the functions $f_{1}, \ldots, f_{s}$
generate the ideal of $Y \setminus V$. Let $M_i \subset
\mathcal{O}(Y)$ be the (finite-dimensional) vector subspace generated
by $G \mal f_i$, and let $N_{i}$ denote the dual $G$-module of
$M_{i}$. Then we obtain $G$-equivariant regular maps
$$ \Phi_i \colon X \to N_i, \qquad x \mapsto (h \mapsto h(x)).$$
Let $N := N_1 \oplus \ldots \oplus N_r$, and let $\Phi \colon Y \to N$
be the map with components $\Phi_i$. Note that $\Phi$ is a
$G$-equivariant closed embedding. Choosing for every $N_i$ a basis
consisting of $H$-homogeneous vectors, we may assume that $N =
\KK^{n}$ holds and that $H$ acts diagonally, i.e., as a subgroup of
the big torus $\TT^{n} \subset \KK^{n}$.

The set $U' \subset \KK^{n}$ consisting of all free $H$-orbits is
invariant under the actions of $G$ and $\TT^{n}$, because these
actions commute with the action of $H$. Moreover,
Remark~\ref{freetorus} implies that $U'$ is open in $\KK^{n}$. Since
$H$ acts freely on $V$, we have $\Phi(V) \subset U'$. Set
$$ U := U' \setminus (N_{s+1} \oplus \ldots \oplus N_r). $$
Then also $U$ is open and invariant under the actions of $G$ and the
big torus $\TT^{n}$. By construction, we have $V = \Phi^{-1}(U)$. So $U$
has the desired properties. The supplement for the case of $G$ being a
torus is obvious. \qed

\medskip

For the next statment, recall that a {\it toric prevariety} is a
normal prevariety together with a regular action of an algebraic torus
having a dense free orbit. An introduction to toric prevarieties is
given in~\cite{acha3}. 

\begin{lemma}\label{torprev}
Notations as in~\ref{cutgeomquots}. The action of $H$ on $U$ admits a
geometric prequotient $p \colon U \to Z := U/H$. Moreover, $Z$ is a
smooth toric prevariety of affine intersection and $G$ acts regularly
on $Z$ making $p  \colon U \to Z$ equivariant.
\end{lemma}

\proof Cover $U$ by $H$-invariant affine open sets $U_{i} \subset U$,
and set $Z_{i} := \Spec(\mathcal{O}(U_{i})^{H})$. Since $H$ acts
freely, the natural maps $p_{i} \colon U_{i} \to Z_{i}$ are geometric
quotients. Using Remark~\ref{freetorus}, one easily verifies that the
maps $p_{i}$ are even locally trivial. In particular, each $Z_{i}$ is
a smooth affine variety.

The varieties $Z_{i}$ glue together along the open subsets
$Z_{ij} := p_{i}(U_{i} \cap U_{j})$ to a smooth prevariety $Z$. Since
each $Z_{ij}$ is the quotient space of the affine variety $U_{i} \cap
U_{j}$, it is again affine. Consequently the prevariety $Z$ is of
affine intersection. Moreover, the maps $p_{i} \colon U_{i} \to Z_{i}$
glue together to a geometric prequotient $p \colon U \to Z$.

Since the actions of $G$ and $\TT^{n}$ on $U$ commute with the
action of $H$, universality of geometric prequotients yields regular
actions of $G$ and $\TT^{n}$ on $Z$ making $p \colon U \to Z$
equivariant. In particular, $Z$ becomes a toric prevariety. \qed 

\medskip

As the $G$-action on the toric prevariety $Z$ in the above lemma is
induced by a linear representation of $G$ on $\KK^{n}$, we call it
{\it linear}. We are now ready for the main result of this section:

\begin{proposition}\label{cone2embed}
Let $G$ be a linear algebraic group, and suppose that the $G$-variety
$X$ has a $G$-equivariant affine conoid. Then $X$ admits a closed
$G$-equivariant embedding into a smooth toric prevariety of affine
intersection on which $G$ acts linearly.
\end{proposition}

\proof Let $\b{X}$ be a $G$-equivariant affine conoid over $X$ and let
$q \colon \rq{X} \to X = \rq{X}/H$ denote the associated geometric quotient.
Lemma~\ref{cutgeomquots} yields a $G \times H$-equivariant embedding $\Phi
\colon \b{X} \to \KK^{n}$ and a $G \times H$-invariant open set $U
\subset \KK^{n}$ with $\Phi^{-1}(U) = \rq{X}$ such that $H$ acts
freely on $U$. As we showed in Lemma~\ref{torprev}, the geometric
prequotient $U \to U/H$ exists and $Z := U/H$ is a smooth toric
prevariety of affine intersection. 

By the universal property of geometric prequotients, the restriction
$\Phi \colon \rq{X} \to U$ induces a regular map $X \to Z$ on the
level of quotients. By construction, this map is equivariant with
respect to the induced linear $G$-action on $Z$. Moreover, by
$H$-closedness of the geometric prequotient $U \to Z$, the map $X \to
Z$ is a closed embedding. \qed

\section{Ample groups of line bundles}\label{amplegroups}

In this section we perform our construction of equivariant affine
conoids. The basic tool is a suitable generalization of ample
line bundles: Instead of a single line bundle, we shall use certain
groups of line bundles. First we make precise what we mean by
a group of line bundles.

Let $X$ be a variety and consider a cover $\mathfrak{U} =
(U_{i})_{i \in I}$ of $X$ by open subsets. This cover gives rise to an 
additive group $\Lambda(\mathfrak{U})$ of line bundles on $X$: For
each cocycle $\xi \in Z^{1}(\mathfrak{U},\mathcal{O}_{X}^{*})$, let
$L_{\xi}$ denote the line bundle obtained by gluing the products
$U_{i} \times \KK$ along the maps
$$ (x,z) \mapsto (x, \xi_{ij}(x)z).$$
The sum $L_{\xi} + L_{\eta}$ of two such line bundles is by definition
the line bundle $L_{\xi\eta} = L_{\eta\xi}$. So the set
$\Lambda(\mathfrak{U})$ consisting of all the bundles $L_{\xi}$ is in
fact an abelian group, isomorphic to
$Z^{1}(\mathfrak{U},\mathcal{O}_{X}^{*})$. When we speak of a group of
line bundles on $X$, we think of a subgroup of some group
$\Lambda(\mathfrak{U})$ as above. 

Now, let $\Lambda$ be a finitely generated free group of line bundles
on $X$. In the sequel, we associate to this group of line bundles a
variety $\rq{X}$ over $X$. For each line bundle $L \in
\Lambda$, let $\mathcal{A}_{L}$ denote its sheaf of sections. We
identify $\mathcal{A}_{0}$ with the structure sheaf $\mathcal{O}_{X}$.

The sections of a line bundle $L_{\xi} \in \Lambda$ over an open set
$U \subset X$ are described by families $f_{i} \in \mathcal{O}_{X}(U
\cap U_{i})$ that are compatible with the gluing cocycle $\xi$. Thus,
for any two sections $f \in \mathcal{A}_{L}(U)$ and $f' \in
\mathcal{A}_{L'}(U)$, we can take the product $(f_{i}f'_{i})$ of their
defining families $(f_{i})$ and $(f'_{i})$ to obtain a
section $ff' \in \mathcal{A}_{L+L'}(U)$. Extending this operation
yields a multiplication on  
$$ \mathcal{A} := \bigoplus_{L \in \Lambda} \mathcal{A}_{L}. $$

We call $\mathcal{A}$ the graded $\mathcal{O}_{X}$-algebra associated
to $\Lambda$. This algebra is reduced, and moreover, it is locally of
finite type over $\mathcal{A}_{0} = \mathcal{O}_{X}$; that means over
sufficently small affine open sets $U \subset X$, the
$\mathcal{O}(U)$-algebra $\mathcal{A}(U)$ is finitely
generated. Consequently we obtain a variety 
$$ \rq{X} := \Spec(\mathcal{A}) $$
by glueing the affine varieties
$\Spec(\mathcal{A}(U))$, where $U$ ranges over small open affine
neighbourhoods $U \subset X$. In this process, the inclusion map
$\mathcal{O}_{X} = \mathcal{A}_{0} \to \mathcal{A}$ gives rise to
an affine regular map
$$ q \colon \rq{X} \to X, $$
and we have
$\mathcal{A} = q_{*}(\mathcal{O}_{\rq{X}})$. We refer to $\rq{X}$ as
to the variety over $X$ associated to the group $\Lambda$. The
$\Lambda$-grading of the $\mathcal{O}_{X}$-algebra $\mathcal{A}$
defines a regular action of the algebraic torus
$$H :=\Spec(\KK[\Lambda])$$
on $\rq{X}$ such that for each affine open set $U \subset
X$, the sections $\mathcal{A}_{L}(U)$ are precisely the functions
of $q^{-1}(U)$ that are homogeneous with respect to the character
$\chi^{L} \in \Chi(H)$. Using Remark~\ref{freetorus}, we observe:

\begin{remark}\label{grothendieck}
\begin{enumerate}
\item $H$ acts freely on $\rq{X}$, and the map $q \colon \rq{X} \to X$
is a geometric quotient for the action of $H$ on $\rq{X}$.
\item For a section $f \in \mathcal{A}_{L}(X)$, let $Z(f) \subset X$
  denote its set of zeroes. Then the set of zeroes of $f$, viewed as a
  regular function on $\rq{X}$, is just
$$ N(\rq{X}; f) \; = \; q^{-1}(Z(f)) \; \subset \; \rq{X}. $$
\end{enumerate}
\end{remark}

To proceed in our construction of affine conoids, we need a condition
on the group $\Lambda$ of line bundles which guarantees that
the associated variety $\rq{X}$ over $X$ is quasiaffine. 

\begin{definition}\label{ampledef}
We call a finitely generated free group $\Lambda$ of line bundles on
$X$ {\it ample} if its associated graded $\mathcal{O}_{X}$-algebra
$\mathcal{A}$ admits homogeneous sections $f_1, \ldots, f_r \in
\mathcal{A}(X)$ such that the open sets $X \setminus Z(f_{i})$ form an
affine cover of $X$.
\end{definition}

Note that this generalizes the usual concept of ampleness in the
sense that an ample line bundle generates an ample group. Moreover, the
above notion of an ample group yields precisely what we are looking for:

\begin{proposition}\label{ample2conoid}
Let $G$ be a linear algebraic group, and let $X$ be a $G$-variety. If
$\Lambda$ is an ample group of line bundles on $X$ and every $L \in
\Lambda$ is $G$-linearizable, then $X$ admits a $G$-equivariant affine
conoid.
\end{proposition}

For the proof we need two statements ensuring existence of 
suitable equivariant affine closures. We use the
following common notation: For a variety $Y$ and a function $h \in
\mathcal{O}(Y)$, let $Y_{h} := \{y \in Y; \; h(y) \ne 0\}$.

\begin{lemma}\label{equivhull}
Let $Y$ be a variety endowed with an action of a linear
algebraic group $G$. Suppose that
\begin{enumerate}
\item there are $f_{1}, \ldots, f_{r} \in \mathcal{O}(Y)$ such
  that the sets $Y_{i} := Y_{f_{i}}$ are affine, cover $Y$ and
  satisfy $\mathcal{O}(Y_{i}) = \mathcal{O}(Y)_{f_{i}}$,
\item the representation of $G$ on $\mathcal{O}(Y)$ given by
$(g \mal f)(y) = f(g^{-1} \mal y)$ is rational.
\end{enumerate}
Then there is an affine $G$-variety $\b{Y}$ containing $Y$ as a dense
open invariant subvariety such that the $f_{i}$ extend regularly to
$\b{Y}$, and $\b{Y}_{f_{i}} = Y_{f_{i}}$ holds. 
\end{lemma}

\proof The main point is that $\mathcal{O}(Y)$ needs not be of
finite type over $\KK$. However, we find $h_{1}, \ldots, h_{s} \in
\mathcal{O}(Y)$ such that there are generators for each $\KK$-algebra
$\mathcal{O}(Y_i)$ among the functions $h_{j}/f_{i}^{l}$. Consider the
subalgebra $A \subset \mathcal{O}(Y)$ generated by $h_{1}, \ldots,
h_{s}$ and $f_{1}, \ldots, f_{r}$. By rationality of the
$G$-representation on $\mathcal{O}(Y)$, we can enlarge $A$ such that
it is $G$-invariant but remains finitely generated.

Consider the affine $G$-variety $\b{Y} := \Spec(A)$. Then the
inclusion $A \subset \mathcal{O}(Y)$ defines a $G$-equivariant regular
map $\varphi \colon Y \to \b{Y}$. Moreover, each function $f_{i} \in
\mathcal{O}(Y)$ is the pullback of a function on $\b{Y}$, denoted
again by $f_{i}$. By construction, restricting
$\varphi$ gives isomorphisms $Y_{f_{i}} \to \b{Y}_{f_{i}}$ of affine
varieties. Since $Y_{f_{i}} = \varphi^{-1}(\b{Y}_{f_{i}})$ holds, we see
that $\varphi$ is the desired open embedding. \qed

\begin{lemma}\label{localizing}
Let $\Lambda$ be an ample group of line bundles on a variety $X$, and
suppose that the sections $f_{1}, \ldots, f_{r}$ of the graded
$\mathcal{O}_{X}$-algebra $\mathcal{A}$ associated to $\Lambda$ are as
in~\ref{ampledef}. Then, setting $X_{i} := X \setminus Z(f_{i})$, we
have $\mathcal{A}(X_{i}) = \mathcal{A}(X)_{f_{i}}$.
\end{lemma}

\proof Let $L_{i} \in \Lambda$ be the degree of $f_{i}$. Then
there is an inverse $f_{i}^{-1} \in \mathcal{A}_{-L_{i}}(X_{i})$ of
$f_{i} \vert_{X_{i}}$. So we obtain an injection
$\mathcal{A}(X)_{f_{i}} \subset \mathcal{A}(X_{i})$. This map is also
surjective: Let $L \in \Lambda$ and $f \in
\mathcal{A}_{L}(X_{i})$. Arguing locally, we see that for some
suitably large integer $m$, the section $f f_{i}^{m} \in
\mathcal{A}_{L + mL_{i}}(X_{i})$ admits an extension to a section of
$\mathcal{A}_{L + mL_{i}}(X)$, compare e.g.~\cite[Proposition
2.2]{Bo}. But this means $f \in \mathcal{A}(X)_{f_{i}}$. \qed

\medskip

\proof[Proof of Proposition~\ref{ample2conoid}]
Fix a basis $L_{1}, \ldots, L_{k}$ of the free abelian group $\Lambda$,
and choose for every $L_{j}$ a $G$-linearization.
Via tensoring these $G$-linearizations, we obtain
a $G$-li\-near\-iza\-tion for each $L \in \Lambda$. This makes the
associated graded $\mathcal{O}_{X}$-algebra $\mathcal{A}$ into a
$G$-sheaf: for a section $f \in \mathcal{A}_{L}(U)$ let
$$ (g \mal f) (x) := g \mal (f(g^{-1} \mal x)).$$
Then $g \mal f \in \mathcal{A}_{L}(g \mal U)$. Moreover, on
$\mathcal{A}_{0} = \mathcal{O}_{X}$ we have the canonical $G$-sheaf
structure arising from the $G$-action on $X$, and the
multiplication of the graded $\mathcal{O}_{X}$-algebra $\mathcal{A}$
associated to $\Lambda$ is compatible with the $G$-action. Since $G$
respects the homogeneous components $\mathcal{A}_{L}$, we infer from
\cite[Section 2.5]{DMV} that the representation of $G$ on
$\mathcal{A}(X)$ is rational. 

Let $\rq{X} = \Spec(\mathcal{A})$ be the variety over $X$ associated
to $\Lambda$, and denote by $q \colon \rq{X} \to X$ the geometric
quotient for the action of $H := \Spec(\KK[\Lambda])$. The fact that
we made $\mathcal{A}$ into a $G$-sheaf, allows us to define a
$G$-action on $\rq{X}$: Fix $g \in G$ and let $U \subset X$ be an
affine open set. The algebra homomorphism 
$$ \mathcal{A}(g \mal U) \to \mathcal{A}(U), \qquad %
f \mapsto g^{-1} \mal f$$ 
defines a morphism of affine varieties
$$ T_{g,U} \colon q^{-1}(U) = \Spec(\mathcal{A}(U)) %
\to \Spec(\mathcal{A}(g \mal U)) = q^{-1}(g \mal U). $$
The maps $T_{g,U}$ glue together to a map $T_{g}
\colon \rq{X} \to \rq{X}$. Moreover, $g \mal x := T_{g}(x)$
defines a group action on $\rq{X}$ and the representation of $G$
on $\mathcal{O}(\rq{X}) = \mathcal{A}(X)$ induced by this action is
the one we started with. In particular, since $\mathcal{A}_{0}$ is
canonically $G$-linearized, the map $q \colon \rq{X} \to X$ is
$G$-equivariant.

In order to check that the actions of $G$ and $H$
on $\rq{X}$ commute, let $x \in \rq{X}$, $g \in G$ and $t \in H$.
Choose a $q$-saturated affine open neighbourhood $\rq{U}$ of $x$. By
$G$-equivariance of $q$, both points $g \mal t \mal x$ and $t \mal g
\mal x$ lie in $g \mal \rq{U}$. Suppose $f \in \mathcal{O}(g \mal
\rq{U})$ is homogeneous with respect to a character $\chi^{L} \in
\Chi(H)$. Then also $g^{-1} \mal f \in \mathcal{O}(\rq{U})$ is
$\chi^{L}$-homogeneous, and we obtain
$$ f(g \mal t \mal x) = (g^{-1} \mal f)(t \mal x)  = \chi^{L}(t)(g^{-1}
\mal f)(x) = \chi^{L}(t)f(g \mal x) = f(t \mal g \mal x). $$
Since the $H$-homogeneous functions separate the points of $g \mal
\rq{U}$, it follows that $g \mal t \mal x$ equals $t \mal g \mal
x$. So the actions of $G$ and $H$ on $\rq{X}$ commute. In particular,
they define an action of the product $G \times H$ on $\rq{X}$.

We shall apply Lemma~\ref{equivhull} to obtain a $G\times
H$-equivariant affine closure $\b{X}$ of $\rq{X}$. First note that
the representation of $G \times H$ on $\mathcal{A}(X)$ is 
rational, because this holds for the representations of the factors $G$
and $H$. So we only have to check Condition~\ref{equivhull}~i).

As to this, let $f_{1}, \ldots, f_{r} \in \mathcal{A}(X)$ as
in~\ref{ampledef}. By Remark~\ref{grothendieck}~ii), it suffices to
know that for $X_{i} := X \setminus Z(f_{i})$, one has
$\mathcal{A}(X_{i}) = \mathcal{A}(X)_{f_{i}}$. But this is guaranteed by
Lemma~\ref{localizing}. So Condition~\ref{equivhull}~i) is verified, and
Lemma~\ref{equivhull} provides a $G \times H$-equivariant affine
closure $\b{X}$ of $\rq{X}$, which is the desired affine conoid. \qed

\medskip

As mentioned earlier, our approach generalizes known constructions for
toric varieties. The most recent one is due to T.~Kajiwara~\cite{Ka};
he proved that a toric variety with enough invariant 
effective Cartier divisors is a geometric quotient of a quasiaffine
toric variety. A systematic treatment of quotient presentations of 
toric varieties is given in~\cite{achass}.

The following observation is useful to construct ample groups of line
bundles on smooth and, more generally $\QQ$-factorial varieties, i.e.,
normal varieties such that for every Weil divisor some multiple is
Cartier. 

\begin{remark}
Suppose on a variety $X$ exist effective Cartier divisors $D_{1},
\ldots ,D_{r}$ such that the sets $X \setminus \Supp(D)_{i}$
form an affine cover of $X$. Then any choice of local equations of the
$D_{i}$ defined on a common cover of $X$ associates to each $D_{i}$ a
line bundle $L_{i}$. Replacing the $D_{i}$ with suitable multiples, 
one achieves that the $L_{i}$ generate a free group $\Lambda$. This
group $\Lambda$ is ample and the canonical sections $f_{i} := 1 \in
\mathcal{O}_{X}(D_{i})$ satisfy to the conditions~\ref{ampledef}.
\end{remark}

\section{Characterizing existence of affine conoids}\label{charconoids}

Summing up the considerations of the preceding two sections, we
present here our first main theorems. Moreover, we give some
discussion and outline a consequence to quotient constructions in this
section. 

Following Borelli~\cite{Bo}, we call a prevariety $X$ {\it divisorial}
if for every $x \in X$ there exists a line bundle $L$ on $X$ admitting
a global section $f$ such that removing its zero set $Z(f)$ yields an
affine open neighbourhood $X \setminus Z(f)$ of $x$.

\begin{remark}
\begin{enumerate}
\item A variety is divisorial if and only if it admits an ample group
  of line bundles.
\item An irreducible variety $X$ is divisorial if and only if for every
  $x \in X$ there exists an effective Cartier divisor $D$ on $X$ such
  that $X \setminus \Supp(D)$ is an open affine neighbourhood of $x$.
\end{enumerate}
\end{remark}

The first main result relates divisoriality, existence of affine
conoids and embeddability into smooth toric prevarieties to each
other:

\begin{theorem}\label{goal1}
For a variety $X$, the following statements are equivalent:
\begin{enumerate}
\item $X$ is divisorial.
\item There exists an affine conoid over $X$.
\item $X$ admits a closed embedding into a smooth toric prevariety of
  affine intersection.
\end{enumerate}
\end{theorem}

\proof The implication ``i)$\Rightarrow$ii)'' is
Proposition~\ref{ample2conoid}, and the implication
``ii)$\Rightarrow$iii)'' is Proposition~\ref{cone2embed}. To obtain
the remaining direction ``iii)$\Rightarrow$i)'', note that a smooth
prevariety of affine intersection is divisorial and that subvarieties
of divisorial prevarieties are again divisorial. \qed 

\begin{remark}
\begin{enumerate}
\item Theorem~\ref{goal1} holds as well for prevarieties $X$. Our
  proof works without changes. Note that a divisorial prevariety $X$
  is necessarily of affine intersection. 
\item The Hironaka twist is a smooth variety of dimension three that
  cannot be embedded into a separated toric variety, compare~\cite{Wl}.
\item There exist normal surfaces that admit neither embeddings into
  toric prevarieties of affine intersection nor into $\QQ$-factorial
  ones, see~\cite{hasc}. 
\end{enumerate}
\end{remark}

\begin{theorem}\label{wlodeq}
Let $X$ be a normal divisorial variety with a regular action
of a connected linear algebraic group $G$. Then there exist a smooth
toric prevariety $Z$ of affine intersection with a linear $G$-action and
a $G$-equivariant closed embedding $X \to Z$.
\end{theorem}

\proof Choose an ample group $\Lambda$ of line bundles on $X$ and fix
a basis $L_{1}, \ldots, L_{r}$ of $\Lambda$. Replacing every $L_{i}$
by a suitable multiple, we achieve that $L_{1}, \ldots, L_{r}$
are $G$-linearizable, see e.g.~\cite[Proposition~2.4]{DMV}. Since
$\Lambda$ remains ample, the assertion follows from
Propositions~\ref{ample2conoid} and~\ref{cone2embed}. \qed

\begin{remark}\label{equembprops}
\begin{enumerate}
\item A condition on $X$ like normality is necessary in~\ref{wlodeq}:
  Identifying $0$ and $\infty$ in the projective line yields
  a $\KK^{*}$-variety that cannot be equivariantly embedded into any
  normal prevariety.
\item If in the setting of~\ref{wlodeq}, the group $G$ is a torus and
acts effectively, then, by the supplement of~\ref{cutgeomquots}, one
can arrange the embedding in such a way that $G$ acts as a subtorus of
the big torus of $Z$. 
\item For $G=\CC^*$, the main result of~\cite{Ha} provides existence
  of equivariant embeddings into toric prevarieties even for
  non-divisorial normal $X$.
\end{enumerate}
\end{remark}

In the remainder of this section we discuss further aspects of affine
conoids. First we note that the results hold also with a more general
definition: In~\ref{conoiddef} we could replace the algebraic torus $H$ by
an arbitrary diagonalizable group. Moreover, in characteristic zero
one could even omit in~\ref{conoiddef} the requirement of $H$ acting
freely on the set $\rq{X}$. This is due to the following observation:

\begin{remark}
Let $H$ be a diagonalizable group acting regularly and effectively on
an affine variety $\b{Y}$. Suppose $\rq{Y} \subset \b{Y}$ is an open
subset with geometric quotient $\rq{Y} \to X := \rq{Y}/H$. Then the
group $\Gamma \subset H$ generated by $H_{y}$, $y \in \rq{Y}$, is
finite. If ${\rm char}(\KK) = 0$, then $\b{X} := \b{Y}/\Gamma$ is an affine
conoid over $X$.
\end{remark}

If a complete variety $X$ admits an ample divisor $D$ in the classical
sense, then the linear system of a suitable multiple of $D$ gives rise
to an embedding of $X$ into a projective space. In the following
example we discuss the embedding of $X$ provided by the (ample) group
$\Lambda$ of line bundles induced by $D$:  

\begin{example}
Let $X$ be a variety with $\mathcal{O}(X) = \KK$, e.g. a complete one.
Suppose that there is a line bundle $L$ on $X$ generating an ample
group $\Lambda = \ZZ L$. We show that the method
of~\ref{ample2conoid} and~\ref{cone2embed} embeds $X$ into a smooth
quasiprojective toric variety:

Let $\rq{X}$ denote the variety over $X$ associated to
$\Lambda$. Note that $\Spec(\KK[\Lambda])$ equals $\KK^{*}$.
Choose any $\KK^{*}$-equivariant affine closure $\b{X}$ of
$\rq{X}$. Since $\mathcal{O}(X) = \KK$ holds, every
$\KK^{*}$-invariant regular function on $\b{X}$ is constant. In
particular, the $\KK^{*}$-variety $\b{X}$ has an attractive fixed
point.

It follows that the map $\Phi$ constructed in Lemma~\ref{cutgeomquots}
embeds $\b{X}$ into some $\KK^{n}$ with linear $\KK^{*}$-action
having zero as attractive fixed point. Hence the induced map $X \to Z$
used in the proof of~\ref{cone2embed} embeds $X$ into the set of
regular points of a weighted projective space. In particular, $X$ is
quasiprojective.
\end{example}

In view of this observation, it is interisting to know when there
exist ``small'' affine conoids over a given variety $X$. For this one
needs small ample groups. Here the Picard group $\Pic(X)$ gives some
bound:

\begin{proposition}
Let $X$ be a divisorial variety. If $\Pic(X)$ is generated by $d$
elements, then $X$ admits an affine conoid $\b{X}$ with $\dim(\b{X})
\le \dim(X) + d$.
\end{proposition}

\proof Choose an ample group $\Lambda$ of line bundles on $X$. Since
$\Pic(X)$ is generated by $d$ elements, there is a subgroup $\Lambda'
\subset \Lambda$ of rank at most $d$ such that each $L \in \Lambda$ is
isomorphic to an element of $\Lambda'$. The variety $\rq{X}$ over $X$
associated to $\Lambda'$ satisfies
$$ \dim(\rq{X}) = \dim(X) + \dim(\Spec(\KK[\Lambda])) = \dim(X) + {\rm
  rk}(\Lambda'). $$
Now, the group $\Lambda'$ is obviously ample. Consequently $\rq{X}$ is
quasiaffine, and any equivariant affine closure $\b{X}$ of $\rq{X}$ is
an affine conoid as wanted. \qed

\medskip

We conclude this section with a ``philosophical'' consequence of
existence of equivariant affine conoids. Assume that a
reductive group $G$ acts regularly on a normal divisorial variety
$X$. It is the central task of Geometric Invariant Theory to look for
$G$-invariant open subsets $U \subset X$ admitting reasonable
quotients. Affine conoids reduce this problem to the quasiaffine case:

\begin{remark}\label{GIT}
Let $\b{X}$ be a $G$-equivariant affine conoid over $X$ and let $q
\colon \rq{X} \to X = \rq{X}/H$ be the associated geometric quotient.
A $G$-invariant open subset $U \subset X$ admits a categorical (good,
geometric) quotient for the action of $G$ if and only if $q^{-1}(U)
\subset \rq{X}$  admits categorical (good, geometric) quotient for the
action of $G \times H$.
\end{remark}

\section{A finiteness result}\label{akproperty}

So far we characterized divisoriality of a given variety $X$ by
existence of an embedding into a smooth toric prevariety $Z$ of affine
intersection. In this section we provide an important ingredient for
the investigation of embeddings into a separated ambient space $Z$.

The following property, also considered in~\cite{Wl} and~\cite{Sw}, is
crucial: We say that a prevariety $X$ has the {\it $A_{k}$-property},
if any $k$ points $x_{1}, \ldots, x_{k} \in X$ admit a common open
affine neighbourhood in $X$. 

\begin{remark}
\begin{enumerate}
\item For $k \ge 2$, an $A_{k}$-prevariety is necessarily
separated.
\item A toric prevariety is separated if and only if it has the
  $A_{2}$-property.
\end{enumerate}
\end{remark}

We are interested in open {\it $A_{k}$-subsets} of a given
prevariety $X$, i.e., open subsets $X' \subset X$ that have as a
prevariety themselves the $A_{k}$-property. The main result of this
section generalizes~\cite[Theorem~3.5]{Sw} to the nonseparated case:

\begin{proposition}\label{Akfinite}
A prevariety has only finitely many maximal open
$A_{k}$-subsets.
\end{proposition}

This proposition can be obtained by combining~\cite[Theorem~3.5]{Sw}
with~\cite[Theorem~I]{BB}. However, for the sake of
self-containedness, we present below a simple direct proof, based on a
slight modification of the arguments used in~\cite[Section~3]{Sw}.  

We apply Proposition~\ref{Akfinite} to actions of connected algebraic
groups $G$. Assume that $G$ acts by means of a regular map $G \times X
\to X$ on a prevariety $X$. As immediate consequences of
the above result, we obtain:

\begin{corollary}\label{maxakinvar}
The maximal open $A_{k}$-subsets of $X$ are $G$-invariant.
\end{corollary}

\proof Compare~\cite[Proof of Prop. 1.3]{acha3}. Let $X_{1}, \ldots,
X_{r}$ be the maximal open $A_{k}$-subsets of $X$. We show that
$X_{1}$ is $G$-invariant. Each $g \in G$ permutes the complements
$A_{i} := X \setminus X_{i}$. Consequently $G$ is covered by the
closed subsets
$$G(i) := \{g \in G; \; g \mal A_{1} \subset A_{i}\}. $$
Now, $G$ is connected, hence $G = G(i)$ for some $i$. In particular,
for the neutral element $e_{G} \in G$ we have $e_{G} \mal A_{1}
\subset A_{i}$. This means $A_{i} = A_{1}$. In other words, $G$ leaves
$X_{1}$ invariant. \qed

\begin{corollary}\label{Gvarak}
$X$ has the $A_{k}$-property if and only if for any collection $B_{1},
\ldots, B_{k} \subset X$ of closed $G$-orbits there exist $x_{i} \in
B_{i}$ such that $x_{1}, \ldots, x_{k}$ admit a common open affine
neighbourhood in $X$. \qed
\end{corollary}

We turn to the proof of Proposition~\ref{Akfinite}. Fix an integer
$k$. Suppose that $X$ is a topological space such that 
the product topology on $X^{k}$ is noetherian, e.g., $X$ is a
prevariety. Let $\mathfrak{U}$ be any family of open subsets of
$X$. Set
$$ A := X^{k} \setminus \bigcup_{U \in \mathfrak{U}} U^{k}. $$
Then $A$ is a closed subspace of $X^{k}$. Denote by $A_{1}, \ldots,
A_{r}$ the irreducible components of $A$. Let $p_{i} \colon X^{k} \to
X$ be the projection onto the $i$-th factor. For a subset $Y \subset
X$, let
$$ X(Y) := X \setminus \bigcup_{p_{i}(A_{j}) \cap Y = \emptyset}
\b{p_{i}(A_{j})}.$$
By a {\it $\mathfrak{U}_{k}$-subset} we mean a subset $Y \subset X$ 
such that for any $x_{1}, \ldots, x_{k} \in Y$ there exists an $U \in
\mathfrak{U}$ that contains the points $x_{1}, \ldots, x_{k}$. The
basic properties of the above construction are subsumed as follows:

\begin{lemma}\label{Ukconstruct}
\begin{enumerate}
\item $X$ has only finitely many subsets of the form $X(Y)$.
\item If $Y$ is open in $X$ then we have $Y \subset  X(Y)$.
\item If $Y \subset X$ is an open $\mathfrak{U}_{k}$-subset then so is
  $X(Y)$.
\end{enumerate}
\end{lemma}

\proof Only for iii) there is something to show. Suppose that $Y$ is
an open $\mathfrak{U}_{k}$-subset but $X(Y)$ does not have the
$\mathfrak{U}_{k}$-property. Then there exist points $x_{1}, \ldots,
x_{k} \in X(Y)$ that are not contained in a common $U \in
\mathfrak{U}$. So $(x_{1}, \ldots, x_{k})$ lies in $A$ and hence in
some irreducible component $A_{j}$ of $A$. In particular, $x_{i}
\in p_{i}(A_{j})$ holds for all $i$.

By definition of $X(Y)$, the fact $x_{i} \in p_{i}(A_{j})$ implies
$p_{i}(A_{j}) \cap Y \ne \emptyset$ for all $i$. Thus each
$p_{i}^{-1}(Y)$ intersects $A_{j}$. Since $A_{j}$ is irreducible and
$Y$ is open, we obtain that $Y^{k}$ intersects $A_{j}$. Since $Y^{k}$
is covered by the sets $U^{k}$, $U \in \mathfrak{U}$, this is a
contradiction to the definition of $A$. \qed

\medskip

\proof[Proof of Proposition~\ref{Akfinite}] Let $\mathfrak{U}$
denote the family of all open affine subvarieties of $X$. According to
Lemma~\ref{Ukconstruct} it suffices to show that the open
$A_{k}$-subsets of $X$ are just its open
$\mathfrak{U}_{k}$-subsets. Clearly every open $A_{k}$-subset $Y
\subset X$ is $\mathfrak{U}_{k}$. The converse is seen as follows:

Let $Y \subset X$ be an open $\mathfrak{U}_{k}$-subset. For given
$x_{1}, \ldots, x_{k} \in Y$ we have to find an affine open $V \subset
Y$ that contains $x_{1}, \ldots, x_{k}$. By assumption, there is an
open affine $U \subset X$ containing $x_{1}, \ldots, x_{k}$. Choose a
function $f \in \mathcal{O}(U)$ that vanishes along $U \setminus Y$
but at no point $x_{i}$. Then $V := U_{f}$ is the desired affine
neighbourhood in $Y$. \qed

\section{Separated ambient spaces}\label{intosep}

Here we discuss embeddings into separated smooth toric varieties. We
shall also consider ambient spaces with additional properties. In order
to formulate our results, we introduce the following terminology:

\begin{definition}\label{kdivdef}
Let $k$ be a positive integer. We say that a prevariety $X$
is {\it $k$-divisorial}, if for any $k$ points $x_{1}, \ldots, x_{k}
\in X$ there is a line bundle $L$ on $X$ admitting
a global section $f$ such that $X \setminus Z(f)$ is affine and
contains $x_{1}, \ldots, x_{k}$. 
\end{definition}

Of course, $k$-divisoriality is strongly related to the
$A_{k}$-property discussed in the preceeding section. Moreover, we
note:

\begin{remark}
\begin{enumerate}
\item A quasiprojective variety is $k$-divisorial for all $k \in \NN$.
\item An irreducible variety $X$ is $k$-divisorial if and only if for
  every $x_{1}, \ldots, x_{k} \in X$ there is an effective Cartier
  divisor $D$ on $X$ such that $X \setminus \Supp(D)$ is affine and
  contains $x_{1}, \ldots, x_{k}$.
\item A $\QQ$-factorial variety is $k$-divisorial if and only if it has
  the $A_{k}$-property.
\item Every $\QQ$-factorial toric variety is $2$-divisorial.
\item The smooth toric variety discussed in~\cite[Example 3.1]{acha2} is
  not $3$-divisorial. 
\end{enumerate}
\end{remark}

The first result of this section characterizes embedabbility into
$k$-divisorial smooth toric varieties. In particular, it implies
\cite[Conjecture~5.3]{Wl} for $\QQ$-factorial varieties:

\begin{theorem}\label{kdivemb}
Let $X$ be a variety, and let $k \ge 2$ be an integer. Then
the following statments are equivalent:
\begin{enumerate}
\item $X$ is $k$-divisorial.
\item $X$ admits a closed embedding into a $k$-divisorial smooth toric
  variety.
\end{enumerate}
\end{theorem}

As a direct consequence, we obtain the following characterization of
embeddability into smooth toric varieties and thereby answer
\cite[Problem~5.4]{Wl} and, partially,~\cite[Problem~5.5]{Wl}:

\begin{corollary}\label{intosepchar}
A variety $X$ admits a closed embedding into a smooth toric
variety if and only if $X$ is $2$-divisorial. \qed
\end{corollary}

As an immediate consequence of this statement, we obtain the following
special version of Nagata's Completion Theorem:

\begin{corollary}
Every $2$-divisorial variety admits a $2$-divisorial completion.
\end{corollary}

\proof Given a $2$-divisorial variety $X$, embed it into a smooth
toric variety $Z$. Choose a smooth toric completion $\b{Z}$ of
$Z$. Then the closure of $X$ in $\b{Z}$ is the desired
completion. \qed

\medskip

Let us turn to $G$-varieties $X$. Though it might be surprising
at the first glance, $k$-divisoriality turns out to be also in the
equivariant setting the right criterion. We prove:

\begin{theorem}\label{kdiveqemb}
Let $G$ be a connected linear algebraic group, and let $X$ be a normal
$G$-variety. If $X$ is $k$-divisorial for some $k \ge 2$, then $X$
admits a $G$-equivariant closed embedding into a smooth $k$-divisorial
toric variety with linear $G$-action.
\end{theorem}

\begin{corollary}\label{V2eq}
Let $G$ be a connected linear algebraic group and let $X$ be
a normal $2$-divisorial $G$-variety. Then $X$ admits a closed
$G$-equivariant embedding of $X$ into a smooth toric variety with
linear $G$-action. \qed
\end{corollary}

As in Remark~\ref{equembprops}~ii), the supplement
of~\ref{cutgeomquots} and the proof given below yield for an effective
action of a torus $G$ on $X$ that one can arrange in
Theorem~\ref{kdiveqemb} the action of $G$ on the ambient toric variety
to be a subtorus action. This implies in particular: 

\begin{corollary}\label{toricemb}
Every $\QQ$-factorial toric variety can be embedded by means of a
toric morphism into a smooth toric variety. \qed
\end{corollary}

\proof[Proof of Theorems~\ref{kdivemb} and~\ref{kdiveqemb}]
By pulling back the desired data from the ambient space, we see that
$k$-divisoriality is necessary to embed a given variety $X$ into
a smooth $k$-divisorial toric variety. We shall show the converse in
the setting of Theorem~\ref{kdiveqemb}. However, normality of $X$ is
merely needed to obtain $G$-linearizations of line bundles. Thus our
proof also settles Theorem~\ref{kdivemb}.

So suppose the connected linear algebraic group $G$ acts regularly
on the normal $k$-divisorial variety $X$. Consider the
$k$-fold product $X^{k}$. This is covered by sets of the form $U^{k}$,
where $U \subset X$ is an affine open subset obtained by removing the
zero set of a section of some line bundle on $X$. Since finitely many
of these $U^{k}$ cover $X^{k}$, we obtain line bundles $L_{1}, \ldots,
L_{r}$ on $X$ and sections $f_{i} \colon X \to  L_{i}$ such that each
$X_{i} := X \setminus Z(f_{i})$ is affine and any $k$ points $x_{1},
\ldots, x_{k} \in X$ lie in some common $X_{i}$. 

Surely, we may assume that $L_{1}, \ldots, L_{r}$ generate a group
$\Lambda$ of line bundles. Moreover, replacing the $L_{i}$ and the
$f_{i}$ with suitable multiples, we achieve that $\Lambda$ is free and
every $L_{i}$ is $G$-linearizable. Then $\Lambda$ is ample, and the sections
$f_{1}, \ldots, f_{r}$ satisfy to the conditions of~\ref{ampledef}. Let
$\rq{X}$ denote the variety over $X$ associated to $\Lambda$. Recall,
that the canonical map $q \colon \rq{X} \to X$ is a geometric quotient
for the action of $H := \Spec(\CC[\Lambda])$ on $\rq{X}$. Moreover, by
Remark~\ref{grothendieck}, $H$ acts freely on $\rq{X}$. 

As in the proof of Proposition~\ref{ample2conoid}, we endow $\rq{X}$
with a $G$-action, commuting with the action of $H$, such that $q
\colon \rq{X} \to X$ becomes $G$-equivariant. In order to obtain
an appropriate $G \times H$-equivariant affine closure
$\b{X}$ of $\rq{X}$, we view the sections $f_{i} \colon X
\to L_{i}$ as regular functions on $\rq{X}$. According
to~\ref{grothendieck} ii), we have $\rq{X}_{f_{i}} = 
q^{-1}(X_{i})$, and Lemma~\ref{localizing} yields
$\mathcal{O}(\rq{X}_{f_{i}}) = \mathcal{O}(\rq{X})_{f_{i}}$. Thus
Lemma~\ref{equivhull} provides a $G \times H$-equivariant affine
closure $\b{X}$ of $\rq{X}$ such that the functions $f_{i}$ extend
regularly to $\b{X}$ and $\b{X}_{f_{i}} = q^{-1}(X_{i})$ holds.
  
Choose a $G \times H$-equivariant embedding $\Phi \colon \b{X} \to
\KK^{n}$ and a $G \times H$-invariant open set $U \subset \KK^{n}$ as
in Lemma~\ref{cutgeomquots}. Then $\Phi^{-1}(U) = \rq{X}$ holds and
the geometric prequotient $U \to Z := U/H$ exists. Moreover, we proved
in Lemma~\ref{torprev} that $Z$ is a smooth toric prevariety
with linear $G$-action. The map $X \to Z$ of quotients induced by
$\Phi$ is a $G$-equivariant closed embedding.

In the sequel, we regard $\b{X}$ and $X$ as subvarieties of $\KK^{n}$
and $Z$ respectively. We claim that for any $k$-points $x_{1}, \ldots,
x_{k} \in X$ there is an affine open neighbourhood $V \subset Z$ with
$x_{1}, \ldots, x_{k} \in V$. To construct such a $V$, we take one of
the $X_{i} \subset X$ with $x_{1}, \ldots, x_{k} \in X_{i}$. By our
choice of $\b{X}$, we have $q^{-1}(X_{i}) = \b{X}_{f_{i}}$.

Now, $f_{i} \in \mathcal{O}(\b{X})$ is the restriction of some
$H$-homogeneous function $h_{i} \in \mathcal{O}(\KK^{n})$.
Consider the $H$-invariant affine open set $U_{i} :=
\KK^{n}_{h_{i}}$. Then $q^{-1}(X_{i})$ is a closed $H$-invariant
subset of $U_{i}$. In particular, $U_{i}$ contains all the fibres
$q^{-1}(x_{j})$. We have to shrink $U_{i}$ a little bit:
Let $A := U_{i} \setminus U$. Then $A$ is a closed $H$-invariant
subset of $U_{i}$. Since $q^{-1}(X_{i}) \subset U$
holds, we obtain $A \cap q^{-1}(X_{i}) = \emptyset$.

Looking at the quotient $\Spec(\mathcal{O}(U_{i})^{H})$, we find a
function $f \in \mathcal{O}(U_{i})^{H}$ that vanishes on $A$ but has
no zeroes along the $H$-orbits $q^{-1}(x_{j})$. Thus,
removing the zero set of this function $f$ from $U_{i}$, we achieve
that $U_{i} \subset U$ holds, $U_{i}$ is still $H$-invariant, affine
and contains all the fibres $q^{-1}(x_{j})$. Now set $V :=
U_{i}/H \subset Z$. Then $V$ is an affine open set in $Z$ and $x_{1},
\ldots, x_{k} \in V$. So our claim is verified.

Let $\SS$ denote the big torus of the toric prevariety $Z$. Removing
from $Z$ step by step the (finitely many) closed $\SS$-orbits that do
not hit $X$, we arrive at an open $\SS$-invariant subset $Z' \subset Z$
such that $X$ is contained in $Z'$ and each closed $\SS$-orbit of $Z'$
has nonempty intersection with $X$. Corollary~\ref{Gvarak} and the
above claim imply that $Z'$ has the $A_{k}$-property and hence is
$k$-divisorial.

To conclude the proof, we have to make $Z'$ invariant under the
action of $G$. We argue in a similar way as above: Let $Z''$ be a
maximal open $A_{k}$-subset of $Z$ such that $Z' \subset Z''$
holds. Since $G$ and $\SS$ are connected, we can apply
Corollary~\ref{maxakinvar}, and obtain that $Z''$ is invariant under
the actions of both, $G$ and $\SS$. So $X \subset Z''$ is the desired
$G$-equivariant closed embedding. \qed

\bibliography{}

\end{document}